\documentclass[11pt,letterpaper,reqno]{amsart}
\usepackage{amssymb,amscd,amsbsy,mathrsfs}
\renewcommand{\a}{\alpha}
\renewcommand{\b}{\beta}
\newcommand{\g}{\gamma}

\newcommand{\e}{\varepsilon}

\newcommand{\s}{\sigma}

\newcommand{\f}{\varphi}

\newcommand{\D}{\Delta}

\newcommand{\E}{{\mathscr E}}

\newcommand{\F}{{\mathscr F}}

\newcommand{\h}{{\mathscr H}}

\newcommand{\K}{{\mathscr K}}

\newcommand{\X}{{\mathscr X}}

\newcommand{\C}{{\Bbb C}}

\newcommand{\R}{{\Bbb R}}
\newcommand{\Z}{{\Bbb Z}}

\newcommand{\bS}{{\boldsymbol S}}

\newcommand{\rf}[1]{(\ref{#1})}

\newcommand{\df}{\stackrel{\mathrm{def}}{=}}

\newcommand{\supp}{\operatorname{supp}}

\newcommand{\trace}{\operatorname{trace}}

\newcommand{\const}{\operatorname{const}}

\newcommand{\eeq}{\end{equation}}
\newcommand{\beq}{\begin{equation}}
\newcommand{\bay}{\begin{eqnarray}}
\newcommand{\ba}{\begin{align*}}
\newcommand{\ea}{\end{align*}}
\newcommand{\ey}{\end{eqnarray}}
\newcommand{\bey}{\begin{eqnarray*}}
\newcommand{\eey}{\end{eqnarray*}}

\newcommand{\be}{\infty}

\newcommand{\bl}{\blacksquare}

\newcommand{\Pf}{{\bf Proof. }}

\newtheorem{thm}{\hspace{\parindent}Theorem}[section]

\newtheorem{cor}[thm]{\hspace{\parindent}Corollary}
\newtheorem{lem}[thm]{\hspace{\parindent}Lemma}

\usepackage{amsmath}
\usepackage{amsfonts}
\usepackage{amssymb}
\usepackage{graphicx}

\usepackage[T2A,T1]{fontenc}
\DeclareSymbolFont{cyrillic}{T2A}{cmr}{m}{it}
\def\makecyrsymbol#1#2{%
    \begingroup\edef\temp{\endgroup
        \noexpand\DeclareMathSymbol{\noexpand#1}
        {\noexpand\mathalpha}{cyrillic}%
        {\expandafter\expandafter\expandafter
            \calccyr\expandafter\meaning\csname T2A\string#2\endcsname\end}}%
    \temp}
\expandafter\def\expandafter\calccyr\string\char#1\end{#1}

\makecyrsymbol\Zhe\CYRZH
\makecyrsymbol\Be\CYRB
\makecyrsymbol\Shcha\CYRSHCH
\makecyrsymbol\Sha\CYRSH

\makeatletter
\def\upintkern@{\mkern-7mu\mathchoice{\mkern-3.5mu}{}{}{}}
\def\upintdots@{\mathchoice{\mkern-4mu\@cdots\mkern-4mu}%
 {{\cdotp}\mkern1.5mu{\cdotp}\mkern1.5mu{\cdotp}}%
 {{\cdotp}\mkern1mu{\cdotp}\mkern1mu{\cdotp}}%
 {{\cdotp}\mkern1mu{\cdotp}\mkern1mu{\cdotp}}}

\newcommand{\UpMultiIntegral}[1]{%
  \edef\ints@c{\noexpand\upintop
    \ifnum#1=\z@\noexpand\upintdots@\else\noexpand\upintkern@\fi
    \ifnum#1>\tw@\noexpand\upintop\noexpand\upintkern@\fi
    \ifnum#1>\thr@@\noexpand\upintop\noexpand\upintkern@\fi
    \noexpand\upintop
    \noexpand\ilimits@
  }%
  \futurelet\@let@token\ints@a
}
\makeatother

\DeclareFontFamily{OMX}{mdbch}{}
\DeclareFontShape{OMX}{mdbch}{m}{n}{ <->s * [0.8]  mdbchr7v }{}
\DeclareFontShape{OMX}{mdbch}{b}{n}{ <->s * [0.8]  mdbchb7v }{}
\DeclareFontShape{OMX}{mdbch}{bx}{n}{<->ssub * mdbch/b/n}{}

\DeclareSymbolFont{uplargesymbols}{OMX}{mdbch}{m}{n}
\SetSymbolFont{uplargesymbols}{bold}{OMX}{mdbch}{b}{n}
\DeclareMathSymbol{\upintop}{\mathop}{uplargesymbols}{82}
\DeclareMathSymbol{\upointop}{\mathop}{uplargesymbols}{"48}

\DeclareFontEncoding{MDB}{}{}
\DeclareFontFamily{MDB}{mdbch}{}
\DeclareFontShape{MDB}{mdbch}{m}{n}{ <->s * [0.8]  mdbchrmb }{}
\DeclareFontShape{MDB}{mdbch}{b}{n}{ <->s * [0.8]  mdbchbmb }{}
\DeclareFontShape{MDB}{mdbch}{bx}{n}{<->ssub * mdbch/b/n}{}
\DeclareFontSubstitution{MDB}{cmr}{m}{n}
\DeclareSymbolFont{mathdesignB}{MDB}{mdbch}{m}{n}%
\SetSymbolFont{mathdesignB}{bold}{MDB}{mdbch}{b}{n}%
\DeclareMathSymbol{\upintclockwise}{\mathop}{mathdesignB}{128}
\DeclareMathSymbol{\upointclockwise}{\mathop}{mathdesignB}{130}
\DeclareMathSymbol{\upointctrclockwise}{\mathop}{mathdesignB}{132}
\DeclareMathSymbol{\upoiint}{\mathop}{mathdesignB}{134}
\DeclareMathSymbol{\upoiiint}{\mathop}{mathdesignB}{136}

\makeatletter
\newcommand{\upint}{\DOTSI\upintop\ilimits@}
\newcommand{\upoint}{\DOTSI\upointop\ilimits@}
\makeatother

\theoremstyle{remark}

\newtheorem*{rem*}{Remark}

\renewcommand{\C}{{\Bbb C}}
\renewcommand{\f}{{\varphi}}
\renewcommand{\b}{{\beta}}

\newcommand\dg{\frak D}

\newcommand\mB{\mathcal{B}}

\newcommand{\ri}{{\rm i}}

\newcommand{\Bs}{\Be_{\be,1}^1}

\newcommand{\CAbe}{{\rm C}_{{\rm A},\be}}

\begin{document}

%
%



\numberwithin{equation}{section}

\numberwithin{equation}{section}

\title{Functons of perturbed pairs\\ of dissipative operators}
\author{A.B. Aleksandrov, V.V. Peller}
\thanks{The research on \S\:4--\S\:6  is supported by 
Russian Science Foundation [grant number 18-11-00053].
The work is supported by a grant of the Government of the Russian Federation for the state support of scientific research, carried out under the supervision of leading scientists, agreement  075-15-2021-602}


\

\begin{abstract}
Let $f$ be a function in the inhomogeneous analytic Besov space $(\Bs)_+(\R^2)$.
For a pair $(L,M)$ of not necessarily commuting maximal dissipative operators,
we define the function
$f(L,M)$ of $L$ and $M$ as a densely defined linear operator.
We prove for $p\in[1,2]$ that if $(L_1,M_1)$ and $(L_2,M_2)$ are pairs of not necessarily commuting maximal dissipative operators
such that both differences $L_1-L_2$ and $M_1-M_2$
belong to the Schatten--von Neumann class $\bS_p$ than for an arbitrary function $f$ in
$(\Bs)_+(\R^2)$, the operator difference
$f(L_1,M_1)-f(L_2,M_2)$ belongs to $\bS_p$ and the following Lipschitz type estimate holds:
$$
\|f(L_1,M_1)-f(L_2,M_2)\|_{\bS_p}
\le\const\|f\|_{\Bs}\max\big\{\|L_1-L_2\|_{\bS_p},\|M_1-M_2\|_{\bS_p}\big\}.
$$
\end{abstract} 

\maketitle


\setcounter{section}{0}
\section{\bf Introduction}
\setcounter{equation}{0}
\label{In}

\

In this paper we obtain Lipschitz type estimates  for functions of pairs of not necessarily commuting maximal dissipative operators 
in the Schatten--von Neumann norm $\bS_p$ with $1\le p\le2$.
Similar results for not necessarily commuting and not necessarily bounded self-adjoint operators were obtained in the paper \cite{AP7}.
The case of commuting dissipative operators was considered by the authors in   \cite{AP6}.

Let $(L,M)$ be a pair of not necessarily commuting maximal dissipative operators. 
Let us define functions  $f(L,M)$ of these operators in the case when the function $f$ on $\R^2$ belongs to the Haagerup tensor product
$\CAbe\!\otimes_{\rm h}\!\CAbe$, where
$$
\CAbe\df H^\be(\C_+)\cap{\rm C}(\R).
$$

Recall that the Haagerup tensor product $\CAbe\!\otimes_{\rm h}\!\CAbe$
consists of functions $f$ of two variables that admit a representation in the form
\bay
\label{Phifnpsin}
f(x,y)=\sum_n\f_n(x)\psi_n(y),
\ey
where $\f_n\in \CAbe$, $\psi_n\in \CAbe$ and
\bay 
\label{koren'izproizvedeniya}
\left(\left\|\sum_n|\f_n(x)|^2\right\|_{L^\be(\R)}
\left\|\sum_n|\psi_n(x)|^2\right\|_{L^\be(\R)}\right)^{1/2}<\be.
\ey
The norm of $f$ in $\CAbe\!\otimes_{\rm h}\CAbe$ 
is, by definition, the infimum of the left-hand side of 
\rf{koren'izproizvedeniya} over all representations of $f$ in the form
\rf{Phifnpsin}.

Let $f\in\!\CAbe\otimes_{\rm h}\!\CAbe$. Then we can define the function  $f$ of the pair 
$(L,M)$ of not necessarily commuting maximal dissipative operators by
\bay
\label{oprfLM}
f(L,M)\df\iint_{\R\times\R}f(x,y)\,d\E_L(x)\,d\E_M(y),
\ey
where $\E_L$ and $\E_M$ are the semi-spectral measures of  $L$ and $M$. We refer the reader to \cite{AP*} and \cite{AP6} where the semi-spectral measure of a maximal dissipative operator is defined and  integrals with respect to semi-spectral measures are discussed.

Note that for $f\in\!\CAbe\otimes_{\rm h}\!\CAbe$, the operator $f(L,M)$ is bounded and
$$
\|f(L,M)\|\le\|f\|_{\CAbe\otimes_{\rm h}\!\CAbe}.
$$

Recall that for a function $f$ in  $\CAbe\otimes_{\rm h}\!\CAbe$, the integral 
on the right-hand side of \rf{oprfLM} can be defined, for example, by
$$
f(L,M)\df\sum_n\f_n(L)\psi_n(M).
$$
The series converges in the weak operator topology and the sum does not depend on the choice of a representation \rf{Phifnpsin}.

In this paper we define the operator $f(L,M)$ for a broader class of functions than  
$\CAbe\otimes_{\rm h}\!\CAbe$, however, in this case the operator $f(L,M)$ is not necessarily bounded anymore. 

Recall that it was shown in \cite{ANP} that for pairs $(A_1,B_1)$ and $(A_2,B_2)$
of not necessarily commuting bounded self-adjoint operators
$A$ and $B$ and a function  $f$  in the 
{\it homogeneous} Besov space $B_{\be,1}^1(\R^2)$ the functions $f(A_1,B_1)$ and $f(A_2,B_2)$ of the operators  $(A_1,B_1)$ and $fA_2,B_2)$
are defined and the following estimate
\bay
\label{vernomezhdu1i2}
\|f(A_1,B_1)-f(A_2,B_2)\|_{\bS_p}\le\const\|f\|_{B_{\be,1}^1}
\max\big\{\|A_1-B_1\|_{\bS_p},\|A_2-B_2\|_{\bS_p}\big\}.
\ey
holds for all $p\in[1,2]$.
In the same paper \cite{ANP} it was shown that for $p>2$ 
(for the operator norm $p=\be$) the inequality \rf{vernomezhdu1i2} cannot hold for all
bounded self-adjoint operators.

The main result of this article is a version of inequality \rf{vernomezhdu1i2} for dissipative operators and for functions $f$ in the {\it inhomogeneous} analytic Besov class $(\Bs)_+(\R^2)$.
In other words, we are going to prove that for $p\in[1,2]$, the following inequality holds:
$$
\|f(L_1,B_1)-f(L_2,M_2)\|_{\bS_p}\le\const\|f\|_{\Bs}
\max\big\{\|L_1-M_1\|_{\bS_p},\|L_2-M_2\|_{\bS_p}\big\}
$$
for all $f\in (\Bs)_+(\R^2)$ and arbitrary dissipative operators $L_1$, $M_1$, $L_2$  and $M_2$.

As in \cite{ANP} an important role here will be played by triple operator integrals. 
A brief introduction ti triple operator integrals is given in \S\,\ref{TrOi}.
In \S\,\ref{nekommu} we define functions of pairs of not necessarily commuting maximal dissipative operators. In \S\,\ref{OtsLiptipa} we obtain for this functional calculus Lipschitz type estimates
 in the norm of $\bS_p$, $1\le p\le2$.

Finally, we give in \S\,\ref{Besovy} a brief discussion of the Besov classes 
$B_{\be,1}^1(\R^2)$,
$\Bs(\R^2)$ and  $(\Bs)_+(\R^2)$.

\

\section{\bf Besov spaces}
\setcounter{equation}{0}
\label{Besovy}

\

In this paper we deal with the homogeneous Besov class $B_{\be,1}^1(\R^2)$, with inhomogeneous Besov class $\Bs(\R^2)$ with the inhomogeneous analytic Besov class 
$(\Bs)_+(\R^2)$.
We refer the reader to the book \cite{Pee}, as well as to the papers  \cite{ANP} and \cite{AP4} for additional information on Besov classes $B_{p,q}^s(\R^d)$. Here we give a definition only in the case when $p=\be$, $q=s=1$ and $d=2$.
Let $w$ be an infinitely differentiable function on the real line $\R$ such that
\bay
\label{w}
w\ge0,\quad\supp w\subset\left[\frac12,2\right],\quad\mbox{and} \quad w(t)=1-w\left(\frac t2\right)\quad\mbox{for}\quad t\in[1,2].
\ey
 
With each $n\in\Z$ we associate the function $W_n$ on $\R^2$ such that
$$
\big(\F W_n\big)(x)=w\left(\frac{\|x\|_2}{2^n}\right),\quad n\in\Z, \quad x=(x_1,x_2),
\quad\|x\|_2\df(x_1^2+x_2^2)^{1/2},
$$
where $\F$ denotes the {\it Fourier transform} defined on $L^1\big(\R^2\big)$ by
$$
\big(\F f\big)(t)=\!\int\limits_{\R^2} f(x)e^{-{\rm i}(x,t)}\,dx,\!\quad 
x=(x_1,x_2),
\quad t=(t_1,t_2), \!\quad(x,t)\df  x_1t_1+x_2t_2.
$$
Clearly,
$$
\sum_{n\in\Z}(\F W_n)(t)=1,\quad t\in\R^2\setminus\{0\}.
$$

With each tempered distribution $f\in{\mathscr S}^\prime\big(\R^2\big)$
we associate the sequence of functions $\{f_n\}_{n\in\Z}$,
\bay
\label{fn}
f_n\df f*W_n.
\ey
The formal series
$
\sum_{n\in\Z}f_n
$
is called a Littlewood--Payley expansion of $f$.
The series does not have to converge to $f$.


First we define the homogeneous Besov class $\dot B^1_{\be,1}\big(\R^2\big)$
as the space of tempered distributions $f\in{\mathscr S}^\prime(\R^2)$ such that
\bay
\label{<be}
\|f\|_{B^1_{\be,1}}\df\sum_{n\in\Z}2^n\|f_n\|_{L^\be}<\be.
\ey
In accordance with this definition the space $\dot B^1_{\be,1}(\R^2)$ contains
all polynomials and $\|f\|_{B^1_{\be,1}}=0$ for every polynomial $f$. Moreover,
a distribution $f$ is determined uniquely by the sequence $\{f_n\}_{n\in\Z}$ up to a polynomial. It is easy to see that the series $\sum_{n\ge0}f_n$ converges in ${\mathscr S}^\prime(\R^2)$.
However, the series $\sum_{n<0}f_n$ in general can diverge. Clearly, the series
\bay
\label{ryad}
\sum_{n<0}\frac{\partial f_n}{\partial x_1}\quad\mbox{and}\quad
\sum_{n<0}\frac{\partial f_n}{\partial x_2}
\ey
converge uniformly on $\R^2$ for an arbitrary $f$ in $\dot B^1_{\be,1}(\R^2)$.
Now we are going to say that a function $f$ belongs to the {\it (reduced) homogeneous Besov class $B^1_{\be,1}(\R^2)$,}
if \rf{<be} holds and
$$
\frac{\partial f}{\partial x_j}=\sum_{n\in\Z}\frac{\partial f_n}{\partial x_j},\quad
j=1,\,2.
$$

In this case the function $f$ is determined by the sequence $\{f_n\}_{n\in\Z}$
uniquely up to a constant, and a polynomial $g$ belongs to $B^1_{\be,1}\big(\R^2\big)$
if and only if it is a constant.

In contrast with our previous papers we are going to use in this paper the inhomogeneous Besov space
$\Bs(\R^2)$. We define the function $W^{[0]}$ by
$$
\big(\F W^{[0]}\big)(t)=1-\sum_{n\ge1}(\F W_n)(t),\quad t\in\R^2. 
$$
It is easy to see that
$$
(\F W^{[0]}\big)(t)=1,\quad\mbox{if}\quad\|t\|_2\le1\quad\mbox{and}\quad
\supp\F W^{[0]}\subset\big\{t\in\R^2:~\|t\|_2\le2\big\}.
$$
We say that a tempered distribution $f$ on $\R^2$
belongs to the {\it inhomogeneous Besov space} $\Bs(\R^2)$ if
\bay
\label{fnol'}
f^{[0]}\df f*W^{[0]}\in L^\be(\R^2)\quad\mbox{and}\quad
\sum_{n\ge1}2^n\|f_n\|_{L^\be}<\be.
\ey
The norm of a function $f$ in $\Bs(\R^2)$ can be defined as
$$
\|f\|_{\Bs}\df\big\|f^{[0]}\big\|_{L^\be}+\sum_{n\ge1}2^n\|f_n\|_{L^\be}.
$$

It is easy to see that if $f\in\Bs(\R^2)$, then
$$
\frac{\partial f}{\partial x_j}=\frac{\partial f^{[0]}}{\partial x_j}
+\sum_{n\ge1}\frac{\partial f_n}{\partial x_j},\quad
j=1,\,2.
$$
Moreover, the series converges uniformly. Clearly,  $\Bs(\R^2)$
is dense in $B^1_{\be,1}(\R^2)$.

We define now the inhomogeneous {\it analytic} Besov class $(\Bs)_+(\R^2)$ by
$$
(\Bs)_+(\R^2)\df\{f\in \Bs(\R^2): \supp\F f\subset[0,+\be)^2\}.
$$
The analytic Besov class $(\Bs)_+(\R^2)$ is a subspace of the Besov class $\Bs(\R^2)$.

We would like to mention that our notation slightly differs form the standard one
(see e.g., \cite{Pee}).
However, we believe that our notation is more convenient.

It is well known and it is easy to prove that  $\Bs(\R^2)=B^1_{\be,1}(\R^2)\bigcap L^\be(\R^2)$. 

With each positive number $\s$ we associate two subspaces
$$
\E_\s^\be(\R^2)\df\Big\{f\in L^\be(\R^2):~\supp\F f\subset\{t\in\R^2:~\|t\|_2\le\s\}\Big\}.
$$
and
$$
(\E_\s^\be)_+(\R^2)\df\Big\{f\in \E_\s^\be(\R^2):\supp\F f\subset[0,+\be)^2\Big\}.
$$
Clearly, $(\E_\s^\be)_+(\R^2)\subset\E_\s^\be(\R^2)\subset{\rm C}^\infty(\R)$.

Note that $f_n\in\E_{2^{n+1}}^\be(\R^2)$ for all $n\in\Z$, where $f_n$ is defined in \rf{fn}.
In particular, if $f\in (\Bs)_+(\R^2)$, then $f_n\in(\E_{2^{n+1}}^\be)_+(\R^2)$
and $f^{[0]}\in(\E_2^\be)_+(\R^2)$, where $f^{[0]}$ is defined in \rf{fnol'}.

\

\section{\bf Triple operator integrals and Haagerup-like tensor products}
\setcounter{equation}{0}
\label{TrOi}

\

A triple operator integral -- is an expression of the form
\bay
\label{troinoi}
\int\limits_{\X_1}\int\limits_{\X_2}\int\limits_{\X_3} 
\Psi(x_1,x_2,x_3)\,dE_1(x_1)T\,dE_2(x_2)R\,dE_3(x_3),
\ey
where $\Psi$ is a bounded measurable function on  $\X_1\times\X_2\times\X_3$; 
$E_1$, $E_2$ and $E_3$ are spectral measures on Hilbert space, and 
$T$ and $R$ are bounded linear operators. Such triple operator integrals
can be defined under certain conditions on $\Psi$, $T$ and $R$.

In  \cite{Pe3} triple operator integrals of the form \rf{troinoi} were defined for arbitrary bounded linear operators $T$ and $R$ and for functions $\Psi$ in the integral projective tensor product  $L^\be(E_1)\!\otimes_{\rm i}\!L^\be(E_2)\!\otimes_{\rm i}\!L^\be(E_3)$ of the spaces $L^\be(E_1)$, 
$L^\be(E_2)$ and $L^\be(E_3)$. In this case the following inequality holds:
$$
 \left\|\iiint\Psi\,dE_1T\,dE_2)R\,dE_3\right\|_{\bS_r}\le
 \|\Psi\|_{L^\be\otimes_{\rm i}L^\be\otimes_{\rm i}L^\be}\|T\|_{\bS_p}\|R\|_{\bS_q},\quad
 T\in\bS_p,~R\in\bS_q,
 $$
 $$
 \frac1r=\frac1p+\frac1q\quad\mbox{under the assumption that}\quad\frac1p+\frac1q\le1.
 $$

Then triple operator integrals were defined in \cite{JTT}
 for the Haagerup tensor products
of the spaces $L^\be$.
{\it The Haagerup tensor product} $L^\be(E_1)\!\otimes_{\rm h}\!L^\be(E_2)\!\otimes_{\rm h}\!L^\be(E_3)$
can be defined as the space of functions $\Psi$ of the form
\bay
\label{htr}
\Psi(x_1,x_2,x_3)=\sum_{j,k\ge0}\a_j(x_1)\b_{jk}(x_2)\g_k(x_3),
\ey
where $\a_j$, $\b_{jk}$ and $\g_k$ are measurable functions such that
\bay
\label{ogr}
\{\a_j\}_{j\ge0}\in L_{E_1}^\be(\ell^2), \quad 
\{\b_{jk}\}_{j,k\ge0}\in L_{E_2}^\be({\mathcal B}),\quad\mbox{and}\quad
\{\g_k\}_{k\ge0}\in L_{E_3}^\be(\ell^2),
\ey
where ${\mathcal B}$ denotes the space of matrices that induce bounded linear operators  $\ell^2$ equipped with the operator norm.

The norm  of $\Psi$ in the space $L^\be\!\otimes_{\rm h}\!L^\be\!\otimes_{\rm h}\!L^\be$
is defined as the infimum of the set of numbers of the form
$$
\|\{\a_j\}_{j\ge0}\|_{L^\be(\ell^2)}\|\{\b_{jk}\}_{j,k\ge0}\|_{L^\be({\mB})}
\|\{\g_k\}_{k\ge0}\|_{L^\be(\ell^2)}
$$
over all representations of $\Psi$ in the form \rf{htr} with families
$\{\a_j\}_{j\ge0}$,
$\{\b_{jk}\}_{j,k\ge0}$ adn $\{\g_k\}_{k\ge0}$ satisfying \rf{ogr}.

Let $\Psi\in L^\be\!\otimes_{\rm h}\!L^\be\!\otimes_{\rm h}\!L^\be$ and let
\rf{htr} and \rf{ogr} be satisfied. Then the triple operator integral
\rf{troinoi} is defined by
\begin{align}
\label{htraz}
\iiint\Psi(x_1,x_2,x_3)&\,dE_1(x_1)T\,dE_2(x_2)R\,dE_3(x_3)\nonumber\\[.2cm]
=&
\sum_{j,k\ge0}\left(\int\a_j\,dE_1\right)T\left(\int\b_{jk}\,dE_2\right)
R\left(\int\g_k\,dE_3\right)\nonumber\\[.2cm]
=&\lim_{M,N\to\be}~\sum_{j=0}^N\sum_{k=0}^M
\left(\int\a_j\,dE_1\right)T\left(\int\b_{jk}\,dE_2\right)
R\left(\int\g_k\,dE_3\right).
\end{align}

The limit in \rf{htraz} exists in the weak operator topology and does not depend on the choice of a representation \rf{htr}, it determines a bounded linear operator and
\bay
\label{ogra}
\left\|\iiint\Psi(x_1,x_2,x_3)\,dE_1(x_1)T\,dE_2(x_2)R\,dE_3(x_3)\right\|
\le\|\Psi\|_{L^\be\!\otimes_{\rm h}\!L^\be\!\otimes_{\rm h}\!L^\be}
\|T\|\cdot\|R\|
\ey
(see \cite{ANP} and \cite{AP5}).

It was proved in \cite{ANP} (see also \cite{AP5}) that if
$$
T\in\bS_p,\quad R\in\bS_q,\quad p,\,q\in[2,\be]\quad\mbox{and}\quad
\Psi\in L^\be\!\otimes_{\rm h}\!L^\be\!\otimes_{\rm h}\!L^\be,
$$
then the triple operator integral \rf{troinoi} belongs to $\bS_r$, $1/r=1/p+1/q$, and
\bay
\label{nerdlyaHaagerupa}
\left\|\iiint\!\Psi(x_1,x_2,x_3)\,dE_1(x_1)T\,dE_2(x_2)R\,dE_3(x_3)\right\|_{\bS_r}
\!\!\!\le\|\Psi\|_{L^\be\!\otimes_{\rm h}\!L^\be\!\otimes_{\rm h}\!L^\be}
\|T\|_{\bS_p}\|R\|_{\bS_q}.
\ey
Note that for $p=\be$ we denote by $\|\cdot\|_\be$ the operator norm.

It turns out that Lipschitz type estimates for functions of pairs of noncommuting operators depend on estimates of triple integral operators with integrands in {\it Haagerup-like tensor products of the first and the second kind} of  $L^\be$ spaces. Such tensor products were introduced in \cite{ANP}.

\medskip

{\bf Definition 1.} A function $\Psi$ of three variables is said to belong to {\it the
Haagerup-like tensor product
$L^\be(E_1)\otimes_{\rm h}L^\be(E_2)\otimes^{\rm h}\!L^\be(E_3)$ of the first kind}
if it admits a representation of the form
\bay
\label{yaH}
\Psi(x_1,x_2,x_3)=\sum_{j,k\ge0}\a_j(x_1)\b_{k}(x_2)\g_{jk}(x_3),\quad x_j\in\X_j,
\ey
with $\{\a_j\}_{j\ge0},~\{\b_k\}_{k\ge0}\in L^\be(\ell^2)$ and
$\{\g_{jk}\}_{j,k\ge0}\in L^\be(\mB)$. As usual, 
$$
\|\Psi\|_{L^\be\otimes_{\rm h}\!L^\be\otimes^{\rm h}\!L^\be}
\df\inf\big\|\{\a_j\}_{j\ge0}\big\|_{L^\be(\ell^2)}
\big\|\{\b_k\}_{k\ge0}\big\|_{L^\be(\ell^2)}
\big\|\{\g_{jk}\}_{j,k\ge0}\big\|_{L^\be(\mB)},
$$
where the infimum is taken over all representations of the form \rf{yaH}.

\medskip

We define now triple operator integrals with integrand in the tensor product
$L^\be(E_1)\!\otimes_{\rm h}\!L^\be(E_2)\!\otimes^{\rm h}\!L^\be(E_3)$.

Let $1\le p\le2$. For a function $\Psi$ in $L^\be(E_1)\otimes_{\rm h}\!L^\be(E_2)\otimes^{\rm h}\!L^\be(E_3)$, for a bounded linear operator $R$ and for an operator $T$ of class $\bS_p$, we define the triple operator integral
\bay
\label{WHft}
W=\iint\!\!\upint\Psi(x_1,x_2,x_3)\,dE_1(x_1)T\,dE_2(x_2)R\,dE_3(x_3)
\ey
as the linear functional on the space $\bS_{p'}$, $1/p+1/p'=1$ (on the space of compact operators in the case $p=1$):
\bay
\label{fko}
Q\mapsto
\trace\left(\left(
\iiint
\Psi(x_1,x_2,x_3)\,dE_2(x_2)R\,dE_3(x_3)Q\,dE_1(x_1)
\right)T\right).
\ey

\medskip

Clearly, the triple operator integral in \rf{fko} is well defined because the function
$$
(x_2,x_3,x_1)\mapsto\Psi(x_1,x_2,x_3)
$$ 
belongs to the Haagerup tensor product  $L^\be(E_2)\!\otimes_{\rm h}\!L^\be(E_3)\!\otimes_{\rm h}\!L^\be(E_1)$.
It follows easily from \rf{nerdlyaHaagerupa} that
\bay
\label{WsSp}
\|W\|_{\bS_p}\le\|\Psi\|_{L^\be\otimes_{\rm h}\!L^\be\otimes^{\rm h}\!L^\be}
\|T\|_{\bS_p}\|R\|,\quad1\le p\le2.
\ey

We also need triple operator integrals in the case when 
$T$ is a bounded linear operator and $R\in\bS_p$, $1\le p\le2$.

\medskip

{\bf Definition 2.} A function $\Psi$ is said to belong to the Haagerup-like tensor product
$L^\be(E_1)\!\otimes^{\rm h}\!L^\be(E_2)\!\otimes_{\rm h}\!L^\be(E_3)$ of the second kind if $\Psi$ admits a representation of the form
\bay
\label{preds}
\Psi(x_1,x_2,x_3)=\sum_{j,k\ge0}\a_{jk}(x_1)\b_{j}(x_2)\g_k(x_3),
\ey
where $\{\b_j\}_{j\ge0},~\{\g_k\}_{k\ge0}\in L^\be(\ell^2)$ and  $\{\a_{jk}\}_{j,k\ge0}\in L^\be(\mB)$.
The norm of $\Psi$ in $L^\be\otimes^{\rm h}\!L^\be\otimes_{\rm h}\!L^\be$ is given by 
$$
\|\Psi\|_{L^\be\otimes^{\rm h}\!L^\be\otimes_{\rm h}\!L^\be}
\df\inf\big\|\{\a_j\}_{j\ge0}\big\|_{L^\be(\ell^2)}
\big\|\{\b_k\}_{k\ge0}\big\|_{L^\be(\ell^2)}
\big\|\{\g_{jk}\}_{j,k\ge0}\big\|_{L^\be(\mB)},
$$
where the infimum is taken over all representations of the form {\em\rf{preds}}.

\medskip

Suppose now that $\Psi\in L^\be(E_1)\!\otimes^{\rm h}\!L^\be(E_2)\!\otimes_{\rm h}\!L^\be(E_3)$,
$T$ is a bounded linear operator and $R\in\bS_p$, $1\le p\le2$. Then the continuous linear functional
$$
Q\mapsto
\trace\left(\left(
\iiint\Psi(x_1,x_2,x_3)\,dE_3(x_3)Q\,dE_1(x_1)T\,dE_2(x_2)
\right)R\right)
$$
on $\bS_{p'}$ (on the space of compact operators in the case $p=1$)
determines an operator $W$ of class $\bS_p$, which we call the triple operator integral \bay
\label{WHst}
W=\upint\!\!\!\iint\Psi(x_1,x_2,x_3)\,dE_1(x_1)T\,dE_2(x_2)R\,dE_3(x_3).
\ey

Clearly,
\bay
\label{eshchorazW}
\|W\|_{\bS_p}\le
\|\Psi\|_{L^\be\otimes^{\rm h}\!L^\be\otimes_{\rm h}\!L^\be}
\|T\|\cdot\|R\|_{\bS_p}.
\ey

In \cite{ANP} more general estimates of triple operator integrals were obtained. Later they were generalized in \cite{AP5}.
 
%

Note that in a similar way we can define Haagerup-like tensor products   
$\CAbe\!\otimes_{\rm h}\!\CAbe\!\otimes^{\rm h}\!\CAbe$ and
$\CAbe\!\otimes^{\rm h}\!\CAbe\!\otimes_{\rm h}\!\CAbe$. 

For a continuously differentiable function $f$ on $\R^2$, we define the divided differences
$\dg^{[1]}f$ and $\dg^{[2]}f$ by
$$
\big(\dg^{[1]}f\big)(x_1,x_2,y)\df\frac{f(x_1,y)-f(x_2,y)}{x_1-x_2}
\quad\mbox{if}\quad x_1\ne x_2,
$$
and
$$
\big(\dg^{[2]}f\big)(x,y_1,y_2)\df\frac{f(x,y_1)-f(x,y_2)}{y_1-y_2}\quad\mbox{if}\quad y_1\ne y_2.
$$
In the case when $x_1=x_2$ or $y_1=y_2$, the divided difference should be replaced with the correspnding partial derivative.

We need the followig fact:

\begin{thm}
\label{anp}
Let $f\in(\E_\s^\be)_+(\R^2)$, $\s>0$. Then $\dg^{[1]}f\in\CAbe\!\otimes_{\rm h}\!\CAbe\!\otimes^{\rm h}\!\CAbe$
and $\dg^{[2]}f\in\CAbe\!\otimes^{\rm h}\!\CAbe\!\otimes_{\rm h}\!\CAbe$; moreover, the following estimates hold:
$$
\|\dg^{[1]}f\|_{\CAbe\otimes_{\rm h}\CAbe\!\otimes^{\rm h}\CAbe}\le\const\s\|f\|_{L^\infty(\R^2)}
$$
and
$$
\|\dg^{[2]}f\|_{\CAbe\otimes^{\rm h}\CAbe\otimes_{\rm h}\CAbe}\le\const\s\|f\|_{L^\infty(\R^2)}.
$$
\end{thm}

An analog  of this result for the spaces $\E_\s^\be(\R^2)$ was obtained in \cite{ANP},
see Corollary 6.3. Theorem \ref{anp} can be proved in a similar way.
First of all, it suffices to consider the case when $\s=1$. Secondly, it suffices to prove the theorem only for the function $\dg^{[1]}$. 
Let us dwell on a single distinction in the proofs. In \cite{ANP} the proof was based on the following identity:
\bay
\label{anp0}
\frac{f(x_1,y)-f(x_2,y)}{x_1-x_2}=\sum_{j,k\in\Z}\frac{\sin(x_1-\pi j)}{(x_1-\pi j)}\cdot\frac{\sin(x_2-\pi k)}{(x_2-\pi k)}
\cdot\frac{f(\pi j,y)-f(\pi k,y)}{\pi(j-k)}
\ey
for all $f\in\E_1^\be(\R^2)$. Theorem \ref{anp} can be proved in a similar way if instead of this identity we use the following one:
\bay
\label{anp1}
\frac{f(x_1,y)-f(x_2,y)}{x_1-x_2}=-\sum_{j,k\in\Z}\frac{e^{{\rm i}x_1}-1}{(x_1-2\pi j)}\cdot\frac{e^{{\rm i}x_2}-1}{(x_2-2\pi k)}
\cdot\frac{f(2\pi j,y)-f(2\pi k,y)}{2\pi(j-k)}
\ey
for all $f\in(\E_1^\be)_+(\R^2)$.
Identity \rf{anp0} follows easily from the following well-known identity:
$$
F(z)=\sum_{n\in\Z}F(\pi n)\frac{\sin(z-\pi n)}{z-\pi n}
$$
for all $F\in L^2(\R)$ such that $\supp \F F\subset[-1,1]$.
Identity \rf{anp0} also follow from equally well-known identity: $$
F(z)=\sum_{n\in\Z}F(2\pi n)\frac{e^{\ri z}-1}{\ri(z-2\pi n)}
$$
for all $F\in L^2(\R)$ such that $\supp\F F\subset[0,1]$.

All the above facts on  operator integrals (including double or triple operator integrals)
also remain valid for semi-spectral measures.

Double operator integrals with respect to semi-spectral measures were defined in 
\cite{Pe2}, see also \cite{AP4} (recall that the difference between the definition of a {\it semi-spectral measure} and the definition of a spectral measure is that instead of the requirement that it takes values in the set of orthogonal projections it is only requires that it takes values in the set of nonnegative contractions, see \cite{AP4} for more detailed information).
In particular, the definition of triple operator integrals  \rf{htraz} for $\Psi$ in the Haagerup tensor product $L^\be(E_1)\!\otimes_{\rm h}\!L^\be(E_2)\!\otimes_{\rm h}\!L^\be(E_3)$ and estimate \rf{ogra} also work in the case when instead of spectral measures $E_1$, $E_2$ and $E_3$  one considers semi-spectral measures $\E_1$, $\E_2$ and $\E_3$. The same can be said about Haagerup-like tensor products of the first and the second kind. We also need the following remark.

\medskip

{\bf Remark.}  Theorems 5.1 and 5.2 of \cite{ANP} also remain valid for {\it semi-spectral} measures  $E_1$, $E_2$ and $E_3$.

\medskip

In what follows we will keep this remark in mind when we refer to Theorems 5.1 and 5.2 of 
\cite{ANP}.

As we have mentioned in the introduction, with each maximal dissipative operator  $L$
one can associate its semi-spectral measure $\E_L$,
defined on the $\s$-algebra of Borel subsets of the real line. Herewith the following
equality holds:
$$
f(L)=\int_\R f(t)\,d\E_L(t)
$$
for every $f$ in $\CAbe$.

We use the notation ${\mathscr B}(\h)$ for the space of bounded linear operators on a Hilbert space $\h$.

If $\K$ is a Hilbert space such that $\h\subset\K$ and $E:{\frak B}\to{\mathscr B}(\K)$ is a spectral measure on a measurable space $(\X,{\frak B})$, it is easy to see that the map 
$\E:{\frak B}\to {\mathscr B}(\h)$ defined by
\bay
\label{dil}
\E(\D)=P_\h E(\D)\big|\h,\quad\D\in{\frak B},
\ey
is a semi-spectral measure. Here $P_\h$ denotes the orthogonal projection onto $\h$.

M.A. Naimark proved in \cite{N} that each semi-spectral measure can be obtained in this way, i.e., each semi-spectral measure is a {\it compression} of a spectral measure. A spectral measure $E$ satisfying \rf{dil} is called a {\it spectral dilation of a semi-spectral measure} $\E$.


\

\         

\section{\bf Functions of not necessarily commuting pairs\\ of maximal dissipative operators}
\label{nekommu}

\

Let $L$ and $M$  be maximal dissipative operators.
For a function $f\in\!\CAbe\!\otimes_{\rm h}\!\CAbe$, we can define the function
$f$ of the pair $(L,M)$ of not necessarily commuting maximal dissipative operators by
$$
f(L,M)\df\iint_{\R\times\R}f(x,y)\,d\E_L(x)\,d\E_M(y).
$$
Obviously, in this case $f(L,M)$ is a bounded linear operator and
$$
\|f(L,M)\|\le\|f\|_{\CAbe\otimes_{\rm h}\!\CAbe}.
$$
Consider the function $f_\sharp$ defined by
$f_\sharp(s,t)\df(1-\ri t)^{-1}f(s,t)$.
Assume that $f_\sharp\in\!\CAbe\otimes_{\rm h}\!\CAbe$. 
Then we can define $f(L,M)$ by
$$
f(L,M)\df f_\sharp(L,M)(I-\ri M)=\left(\,\,\,\iint\limits_{\R\times\R}f_\sharp(s,t)\,d\E_L(s)\,d\E_M(t)\right)(I-\ri M).
$$
Note that the operator $f(L,M)$ is defined on $D(M)$, where $D(M)$ denotes {\it the domain of $M$}. In this case the operator $f(L,M)$ does not have to be bounded but the operator $f(L,M)(I-\ri M)^{-1}$ must be bounded.

By Corollary 7.3 of \cite{AP6}, the inclusion $f_\sharp\in(\E_\s^\be)_+(\R^2)$ holds for an arbitrary function $f$ in $(\E_\s^\be)_+(\R^2)$, where $\s>0$.
Thus, the operator $f_\sharp(L,M)$ is bounded for every function $f$ in $(\E_\s^\be)_+(\R^2)$, where $\s>0$ while the operator $f(L,M)$ is not necessarily a bounded operator
defined on $D(M)$.
Moreover,
\bay
\label{1+sigma+}
\|f_\sharp\|_{\CAbe\otimes_{\rm h}\!\CAbe}\le\const(1+\s)\|f\|_{L^\be(\R^2)},
\ey
see Corollary 7.3 of \cite{AP6}.

\medskip

{\bf Remark.} Let $M_0$ be another maximal dissipative operator such that the operator $M_0-M$ is bounded.
Then $f(L,M)(I-\ri M_0)^{-1}$ is a bounded operator for an arbitrary function  $f$ in $(\E_\s^\be)(\R^2)$, where $\s>0$.

\medskip

\Pf It suffices to observe that
\bey
f(L,M)(I-\ri M_0)^{-1}=f_\sharp(L,M)(I-\ri M)(I-\ri M_0)^{-1}\\
=\ri f_\sharp(L,M)(M_0-M)(I-\ri M_0)^{-1}+f_\sharp(L,M). \quad\bl
\eey

\begin{thm} 
\label{0th73}
Let $f\in\big(\Bs\big)_+(\R^2)$. Then $f_\sharp\in \!\CAbe\otimes_{\rm h}\!\CAbe$
and $\|f_\sharp\|_{\CAbe\otimes_{\rm h}\!\CAbe}\le\const\|f\|_{(\Bs)_+}$.
\end{thm}

\Pf Let $f\in(\Bs)_+(\R^2)\subset\Bs(\R^2)$. Then $f=f^{[0]}+\sum\limits_{n=1}^\infty f_n$ (the series converges in  $(\Bs)_+(\R^2)$), 
where $f^{[0]}$ and $f_n$ denote the same as in formulae \rf{fn} and \rf{fnol'}. It remains to apply inequality \rf{1+sigma+}
to the function $f^{[0]}$ and to all functions $f_n$ with $n\ge1$. $\bl$

\medskip

Theorem \ref{0th73} allows us for an arbitrary function $f$ in $\left(\Bs\right)_+(\R^2)$ to define the operator $f(L,M)$  
with domain $D(M)$ for arbitrary maximal dissipative operators $L$ and $M$.

\

\section{\bf Integral formulae and Lipschitz type estimates for the difference of functions of pairs of maximal dissipative operators}
\label{OtsLiptipa}

\

Let $L$ be a maximal dissipative operator. Put $L(\e)\df L(I-\ri\e L)^{-1}$ for $\e>0$.
It is easy to see that $L(\e)$ is a bounded dissipative operator for every positive number $\e$.

In what follows for not necessarily commuting maximal dissipative operators $L_1$ and $L_2$ we will say that $L_1-L_2\in\bS_p$ if the operator $L_1-L_2$ extends by continuity to an operator of class $\bS_p$.

\begin{lem}
\label{0spnorm} Let $L_1$ and $L_2$ be maximal dissipative operators such that $L_1-L_2\in\bS_2$.
Then
$$
\lim_{\e\to0}(L_1(\e)-L_2(\e))=L_1-L_2
$$
in the norm of $\bS_2$.
\end{lem}

The lemma can be proved by analogy with the proof of Lemma 5.5 of  \cite{AP7}.


\begin{thm}
\label{teor1}
Let $f\in(\E^\be_\s)_+(\R^2)$. Suppose that  $L_1$, $L_2$ and $M$ --
are maximal dissipative operators such that
$L_1-L_2\in\bS_2$. Then the following identity holds:
\begin{align*}
f(L_1,M)&-f(L_2,M)\\[.2cm]
&=
\iint\!\!\upint\frac{f(x_1,y)-f(x_2,y)}{x_1-x_2}
\,d\E_{L_1}(x_1)(L_1-L_2)\,d\E_{L_2}(x_2)\,d\E_{M}(y).
\end{align*}
\end{thm}

{\bf Remark.} Strictly speaking, the left-hand side is defined on $D(M)$ while the right-hand side is defined on the whole space. We mean here that the operator on the left
extends by continuity to the whole space and the extension coincides with the right-hand side.

\medskip

\Pf It follows from Theorem \ref{anp} that the divided difference $\dg^{[1]}f$
belongs to the space $\CAbe\!\otimes_{\rm h}\!\CAbe\!\otimes^{\rm h}\!\CAbe$
and so the triple operator integral is well defined.

To prove the desired equality, it suffices to verify that
\begin{align*}
(f(L_1,M)&-f(L_2,M))(I-\ri M)^{-1}\\[.2cm]
&\!\!\!\!\!=
\left(\iint\!\!\upint\frac{f(x_1,y)-f(x_2,y)}{x_1-x_2}
\,d\E_{L_1}(x_1)(L_1-L_2)\,d\E_{L_2}(x_2)\,d\E_{M}(y)\right)(I-\ri M)^{-1}.
\end{align*}

Note that Lemma \ref{0spnorm} of this paper and Theorem 5.1 of \cite{ANP} imply the following equality:
\begin{align*}
\lim_{\e\to0}\left(\iint\!\!\upint\frac{f(x_1,y)-f(x_2,y)}{x_1-x_2}
\,d\E_{L_1}(x_1)(L_1(\e)-L_2(\e))\,d\E_{L_2}(x_2)\,d\E_{M}(y)\right)(I-\ri M)^{-1}\\
=\iint\!\!\upint\frac{f(x_1,y)-f(x_2,y)}{x_1-x_2}
\,d\E_{L_1}(x_1)(L_1-L_2)\,d\E_{L_2}(x_2)(I-\ri M)^{-1}\,d\E_{M}(y)
\end{align*}
in the norm of $\bS_2$.
On the other hand, we have
\begin{multline*}
\lim_{\e\to0}\left(\iint\!\!\upint\!\frac{f(x_1,y)-f(x_2,y)}{x_1-x_2}\!
\,d\E_{L_1}(x_1)(L_1(\e)-L_2(\e))\,d\E_{L_2}(x_2)\,d\E_{M}(y)\right)(I-\ri M)^{-1}\\
=\!\lim_{\e\to0}\!\!\left(\iint\!\!\upint\!\!\frac{f(x_1,y)\!-\!f(x_2,y)}{x_1-x_2}\!\left(\!\frac{x_1}{1\!-\!\ri\e x_1}\!-\!\frac{x_2}{1\!-\!\ri\e x_2}\!\right)\!
d\E_{L_1}(x_1)d\E_{L_2}(x_2)d\E_{M}(y)\!\right)\!\!(I-\ri  M)^{\!-1}\\
=\lim_{\e\to0}\iint\!\!\upint\frac{f(x_1,y)-f(x_2,y)}{(1-\ri \e x_1)(1-\ri \e x_2)(1-\ri y)}
\,d\E_{L_1}(x_1)\,d\E_{L_2}(x_2)\,d\E_{M}(y)\\
=\lim_{\e\to0}\iint\!\!\upint\frac{f(x_1,y)}{(1-\ri \e x_1)(1-\ri \e x_2)(1-\ri y)}
\,d\E_{L_1}(x_1)\,d\E_{L_2}(x_2)\,d\E_{M}(y)\\
-\lim_{\e\to0}\iint\!\!\upint\frac{f(x_2,y)}{(1-\ri \e x_1)(1-\ri \e x_2)(1-\ri y)}
\,d\E_{L_1}(x_1)\,d\E_{L_2}(x_2)\,d\E_{M}\\
=\lim_{\e\to0}\iint \frac{f(x_1,y)}{(1-\ri \e x_1)(1-\ri y)}
\,d\E_{L_1}(x_1)(1-\ri \e L_2)^{-1}\,d\E_{M}(y)\\
-\lim_{\e\to0}(1-\ri \e L_1)^{-1}\iint\frac{f(x_2,y)}{(1-\ri \e x_2)(1-\ri y)}
\,d\E_{L_2}(x_2)\,d\E_{M}(y)
=\lim_{\e\to0} X_\e-\lim_{\e\to0} Y_\e\,,
\end{multline*}
where
$$
X_\e\df\lim_{\e\to0}(1-\ri \e L_1)^{-1}\iint\frac{f(x_1,y)}{1-\ri y}
\,d\E_{L_1}(x_1)(1-\ri \e L_2)^{-1}\,d\E_{M}(y)
$$
and
$$
Y_\e\df
\lim_{\e\to0}(1-\ri \e L_1)^{-1}(1-\ri \e L_2)^{-1}\iint \frac{f(x_2,y)}{1-\ri y}
\,d\E_{L_2}(x_2)\,d\E_{M}(y).
$$

Clearly, 
$$
\lim_{\e\to0} Y_\e=\iint \frac{f(x_2,y)}{1-\ri y}\,d\E_{L_2}(x_2)\,d\E_{M}(y)=f(L_2,M)(I-\ri  M)^{-1}
$$
in the strong operator topology.
To evaluate the limit $\lim\limits_{\e\to0} X_\e$, we observe that by Proposition
3.3\footnote{In \cite{BS3} the authors considered only spectral measures; however, the proof given there also works in the case of semi-spectral measures.}
of \cite{BS3},
\bay
\label{prop33}
\lim_{\e\to0}\iint\frac{f(x_1,y)}{1-\ri y}\,d\E_{L_1}(x_1)(1-\ri \e L_2)^{-1}\,d\E_{M}(y)
=\iint\frac{f(x_1,y)}{1-\ri y}\,d\E_{L_1}(x_1)\,d\E_{M}(y)
\ey
in the strong operator topology. Therefore,
\begin{align*}
\lim_{\e\to0} X_\e=\iint\frac{f(x_1,y)}{1-\ri y}\,d\E_{L_1}(x_1)\,d\E_{M}(y)
=f(L_1,M)(I-\ri M)^{-1}
\end{align*}
in the strong operator topology. Thus, we have proved that
\bey
\lim_{\e\to0}\left(\iint\!\!\upint\frac{f(x_1,y)-f(x_2,y)}{x_1-x_2}
\,d\E_{L_1}(x_1)(L_1(\e)-L_2(\e))\,d\E_{L_2}(x_2)\,d\E_{M}(y)\right)(I-\ri  M)^{-1}\\
=(f(L_1,M)-f(L_2,M))(I-\ri  M)^{-1}
\eey
in the strong operator topology. $\bl$.

\begin{cor} 
\label{cor57}
Let $f\in(\E^\be_\s)_+(\R^2)$ and let $1\le p\le2$. Suppose that $L_1$, $L_2$ and $M$ are maximal dissipative operators such that
$L_2-L_1\in\bS_p$. Then the following inequality holds:
$$
\|f(L_1,M)-f(L_2,M)\|_{\bS_p}\le\const\s\|f\|_{L^\be(\R^2)}\|L_1-L_2\|_{\bS_p}.
$$
\end{cor}

\Pf Applying estimate  \rf{WsSp} and Theorem \ref{anp}, we obtain
\begin{align*}
\|f(L_1,M)-f(L_2,M)\|_{\bS_p}
&\le\big\|\big(\dg^{[1]}f\big)(x_1,x_2,y)\big\|_{\CAbe\!\otimes_{\rm h}\!\CAbe\!\otimes^{\rm h}\!\CAbe}\|L_1-L_2\|_{\bS_p}\\
&\le\const\s\|f\|_{L^\be(\R^2)}\|L_1-L_2\|_{\bS_p}.\qquad\qquad\qquad\bl
\end{align*}

\begin{cor} 
\label{cor58}
Let $f\in\big(\Bs\big)_+(\R^2)$, $1\le p\le2$. Suppose that  $L_1$, $L_2$ and $M$ are maximal dissipative operators such that 
$L_2-L_1\in\bS_p$. Then the following inequality holds:
$$
\|f(L_1,M)-f(L_2,M)\|_{\bS_p}\le\const\|f\|_{\Bs}\|L_1-L_2\|_{\bS_p}.
$$
\end{cor}

\Pf Let $f^{[0]}$ and $f_n$ are the functions defined by \rf{fnol'} and \rf{fn}.
It is easy to see that if $u\in D(M)$, then
$$
f(L_1,M)u-f(L_2,M)u=f^{[0]}(L_1,M)u-f^{[0]}(L_2,M)u+\sum_{n=1}^\infty(f_n(L_1,M)u-f_n(L_2,M)u).
$$
It remains to observe that
$$
\|f^{[0]}(L_1,M)u-f^{[0]}(L_2,M)u\|\le\const\|f^{[0]}\|_{L^\infty}\|L_1-L_2\|\cdot\|u\|
$$
and
$$
\|f_n(L_1,M)u-f_n(L_2,M)u\|\le\const 2^n\|f_n\|_{L^\infty}\|L_1-L_2\|\cdot\|u\|
$$
by Corollary \ref{cor57}. $\bl$

\begin{thm}
\label{teor2}
Let $f\in(\E^\be_\s)_+(\R^2)$. Suppose that $L$, $M_1$ and $M_2$ are maximal dissipative operators  such that
$M_2-M_1\in\bS_2$. Then the following equality holds: 
\begin{align*}
f(L,M_1)&-f(L,M_2)\nonumber\\[.2cm]
&=
\upint\!\!\!\iint\frac{f(x,y_1)-f(x,y_2)}{y_1-y_2}
\,d\E_{L}(x)\,d\E_{M_1}(y_1)(M_1-M_2)\,d\E_{M_2}(y_2).
\end{align*}
\end{thm}

\Pf The proof is similar to the proof of the proof of Theorem \ref{teor1}. It follows from Theorem \ref{anp}
that the right-hand side is well defined.

To prove the desired equality, it suffices to verify that
\begin{align*}
(f(L,M_1)&-f(L,M_2))(I-\ri M_2)^{-1}\nonumber\\[.2cm]
&\!=\!
\left(\upint\!\!\!\iint\frac{f(x,y_1)-f(x,y_2)}{y_1-y_2}
d\E_{L}(x)d\E_{M_1}(y_1)(M_1-M_2)d\E_{M_2}(y_2)\right)(I-{\rm i}M_2)^{-1}.
\end{align*}

Note that Lemma \ref{0spnorm} of this paper and Theorem 5.2 of  \cite{ANP}
imply the following equality:
\begin{align*}
\lim_{\e\to0}\left(\upint\!\!\!\iint\frac{f(x,y_1)-f(x,y_2)}{y_1-y_2}
\,d\E_{L}(x)\,d\E_{M_1}(y_1)(M_1(\e)-M_2(\e))\,d\E_{M_2}(y_2)\right)(I-{\rm i}M_2)^{-1}\\
=\left(\upint\!\!\!\iint\frac{f(x,y_1)-f(x,y_2)}{y_1-y_2}
\,d\E_{L}(x)\,d\E_{M_1}(y_1)(M_1-M_2)\,d\E_{M_2}(y_2)\right)(I-{\rm i}M_2)^{-1}
\end{align*}
in the norm of $\bS_2$. On the other hand, we have
\begin{multline*}
\lim_{\e\to0}\left(\upint\!\!\!\iint\frac{f(x,y_1)-f(x,y_2)}{y_1-y_2}\,d\E_{L}(x)
\,d\E_{M_1}(y_1)(M_1(\e)-M_2(\e))\,d\E_{M_2}(y_2)\right)(I-{\rm i}M_2)^{-1}\\
=\!\lim_{\e\to0}\left(\upint\!\!\!\iint\!\frac{f(x,y_1)\!-\!f(x,y_2)}{y_1-y_2}\!
\left(\!\frac{y_1}{1\!-\!{\rm i}\e y_1}\!-\!\frac{y_2}{1\!-\!{\rm i}\e y_2}\!\right)
\!d\E_{L}(x)d\E_{M_1}(y_1)d\E_{M_2}(y_2)\right)\!(I-{\rm i}M_2)^{\!-1}\\
=\lim_{\e\to0}\upint\!\!\!\iint\frac{f(x,y_1)-f(x,y_2)}{(1-{\rm i}\e y_1)(1-{\rm i}\e y_2)(1-{\rm i}y_2)}
\,d\E_{L}(x)\,d\E_{M_1}(y_1)\,d\E_{M_2}(y_2)\\
=\lim_{\e\to0}\upint\!\!\!\iint\frac{f(x,y_1)}{(1-{\rm i}\e y_1)(1-{\rm i}\e y_2)(1-{\rm i}y_2)}
\,d\E_{L}(x)\,d\E_{M_1}(y_1)\,d\E_{M_2}(y_2)\\
-\lim_{\e\to0}\upint\!\!\!\iint\frac{f(x,y_2)}{(1-{\rm i}\e y_1)(1-{\rm i}\e y_2)(1-{\rm i}y_2)}
\,d\E_{L}(x)\,d\E_{M_1}(y_1)\,d\E_{M_2}(y_2)=X_\e-Y_\e,
\end{multline*}
where
$$
X_\e\df\lim_{\e\to0}\left(\iint\frac{f(x,y_1)}{1-{\rm i}\e y_1}
\,d\E_{L}(x)\,d\E_{M_1}(y_1)\right)(I-{\rm i}\e M_2)^{-1}(I-{\rm i}M_2)^{-1}
$$
and
$$
Y_\e\df\lim_{\e\to0}\left(\iint\frac{f(x,y_2)}{(1-{\rm i}\e y_2)(1-{\rm i}y_2)}
\,d\E_{L}(x)(I-{\rm i}\e M_1)^{-1}\,d\E_{M_2}(y_2)\right).
$$


Clearly,
\begin{multline*}
\lim_{\e\to0} X_\e\\
=\lim_{\e\to0}\left(\iint\frac{f(x,y_1)}{1-{\rm i}y_1}
\,d\E_{L}(x)\,d\E_{M_1}(y_1)\right)(I-{\rm i}\e M_1)^{-1}(I-{\rm i}M_1)(I-{\rm i}M_2)^{-1}(I-{\rm i}\e M_2)^{-1}\\
=\left(\iint\frac{f(x,y_1)}{1-{\rm i}y_1}
\,d\E_{L}(x)\,d\E_{M_1}(y_1)\right)(I-{\rm i}M_1)(I-{\rm i}M_2)^{-1}\\
=f(L,M_1)(I-{\rm i}M_2)^{-1}
\end{multline*}
in the strong operator topology.
To evaluate the limit $\lim\limits_{\e\to0} Y_\e$, we observe that by analogy with the passage to the limit in \rf{prop33}, we have
\begin{align*}
\lim_{\e\to0}\iint\frac{f(x,y_2)}{1-{\rm i} y_2}
\,d\E_{L}(x)(I-{\rm i}\e M_1)^{-1}\,d\E_{M_2}(y_2)
=\iint\frac{f(x,y_2)}{1-{\rm i} y_2}
\,d\E_{L}(x)\,d\E_{M_2}(y_2)
\end{align*}
in the strong operator topology. Consequently,
\begin{align*}
\lim_{\e\to0} Y_\e=\lim_{\e\to0}\left(\iint\frac{f(x,y_2)}{1-{\rm i} y_2}
\,d\E_{L}(x)(I-{\rm i}\e M_1)^{-1}\,d\E_{M_2}(y_2)\right)(I-{\rm i}\e M_2)^{-1}\\
=\iint\frac{f(x,y_2)}{1-{\rm i} y_2}
\,d\E_{L}(x)\,d\E_{M_2}(y_2)=f(L,M_2)(1-{\rm i}M_2)^{-1}.
\end{align*}
in the strong operator topology. Thus, we have proved that

\begin{align*}
\lim_{\e\to0}\left(\upint\!\!\!\iint\frac{f(x,y_1)-f(x,y_2)}{y_1-y_2}
\,d\E_{L}(x)\,d\E_{M_1}(y_1)(M_1(\e)-M_2(\e))\,d\E_{M_2}(y_2)\right)(I-{\rm i}M_2)^{-1}\\
=(f(L,M_1)-f(L,M_2))(I-{\rm i}M_2)^{-1}
\end{align*}
in the strong operator topology. $\bl$.

\begin{cor}
\label{cor510}
Let $f\in(\E^\be_\s)_+(\R^2)$ and let $1\le p\le2$. Suppose that $L$, $M_1$ and $M_2$ are maximal dissipative operators such that
$M_2-M_1\in\bS_p$. Then the following estimate holds:
$$
\|f(L,M_1)-f(L,M_2)\|_{\bS_p}\le\const\s\|f\|_{L^\be(\R^2)}\|M_1-M_2\|_{\bS_p}.
$$
\end{cor}

\Pf Applying estimate \rf{eshchorazW}  and Theorem \ref{anp}, we obtain
\bey
\|f(L,M_1)-f(L,M_2)\|_{\bS_p}
\le\big\|\big(\dg^{[2]}f\big)(x,y_1,y_2)\big\|_{\CAbe\!\otimes_{\rm h}\!\CAbe\!\otimes^{\rm h}\!\CAbe}\|M_1-M_2\|_{\bS_p}\\
\le\const\s\|f\|_{L^\be(\R^2)}\|L_1-L_2\|_{\bS_p}.\quad\bl
\eey

\begin{cor} 
\label{511}
Let $f\in(\Bs)_+(\R^2)$, $1\le p\le2$. Suppose that $L$, $M_1$ and $M_2$ are maximal dissipative operators such that
 $M_2-M_1\in\bS_p$. Then the following inequality holds:
$$
\|f(L,M_1)-f(L,M_2)\|_{\bS_p}\le\const\|f\|_{\Bs}\|M_1-M_2\|_{\bS_p}.
$$
\end{cor}

\Pf Let $f^{[0]}$ and $f_n$ be the functions 
defined by \rf{fnol'} and \rf{fn}.
It is easy to see that if $u\in D(M_1)=D(M_2)$ then
$$
f(L,M_1)u-f(L,M_2)u=f^{[0]}(L, M_1)u-f^{[0]}(L,M_2)u+\sum_{n=1}^\infty(f_n(L,M_1)u-f_n(L,M_2)u).
$$
It remains to observe that
$$
\|f^{[0]}(L,M_1)u-f^{[0]}(L,M_2)u\|\le\const\|f^{[0]}\|_{L^\infty}\|M_1-M_2\|\cdot\|u\|
$$
and
$$
\|f_n(L,M_1)u-f_n(L,M_2)u\|\le\const 2^n\|f_n\|_{L^\infty}\|M_1-M_2\|\cdot\|u\|
$$
by Corollary \ref{cor510}. $\bl$

\begin{thm}
\label{lipschitseva_otsenka_dlya_Besova}
Let $f\in(\Bs)_+(\R^2)$, $1\le p\le2$. Suppose that $L_1$, $L_2$, $M_1$ and $M_2$ are maximal dissipative operators such that
$L_2-L_1\in\bS_p$ and $M_2-M_1\in\bS_p$. Then the following inequality holds:
$$
\|f(L_1,M_1)-f(L_2,M_2)\|_{\bS_p}\le\const\|f\|_{\Bs}
\max\big\{\|L_1-L_2\|_{\bS_p},\|M_1-M_2\|_{\bS_p}\big\}.
$$
\end{thm}

\Pf We have
$$
\|f(L_1,M_1)-f(L_2,M_2)\|_{\bS_p}\le\|f(L_1,M_1)-f(L_2,M_1)\|_{\bS_p}+\|f(L_2,M_1)-f(L_2,M_2)\|_{\bS_p}.
$$
It remains to apply Corollaries \ref{cor58}  and \ref{511}. $\bl$

\begin{thm} 
\label{ab12}
Let $f\in(\Bs)_+(\R^2)$. Suppose that $L_1$, $L_2$, $M_1$ and $M_2$ are maximal dissipative operators such that
$L_1-L_2\in\bS_2$ and $M_1-M_2\in\bS_2$. Then the following equality holds:
\begin{align}
\label{tryokhetazhka}
f(L_1,M_1)&-f(L_2,M_2)\nonumber\\[.2cm]
&=
\iint\!\!\upint\frac{f(x_1,y)-f(x_2,y)}{x_1-x_2}
\,d\E_{L_1}(x_1)(L_1-L_2)\,d\E_{L_2}(x_2)\,d\E_{M_1}(y),\nonumber\\[.2cm]
&+\upint\!\!\!\iint\frac{f(x,y_1)-f(x,y_2)}{y_1-y_2}
\,d\E_{L_2}(x)\,d\E_{M_1}(y_1)(M_1-M_2)\,d\E_{M_2}(y_2).
\end{align}
\end{thm}

Here it is worth to say the same words as in the remark following Theorem 
 \ref{teor1}.

\medskip

\Pf Suppose first that $f\in \E^\be_\s(\R^2)$.By Theorem \ref{teor1},
we have
\begin{align*}
f(L_1,M_1)&-f(L_2,M_1)\\[.2cm]
&=
\iint\!\!\upint\frac{f(x_1,y)-f(x_2,y)}{x_1-x_2}
\,d\E_{L_1}(x_1)(L_1-L_2)\,d\E_{L_2}(x_2)\,d\E_{M_1}(y),
\end{align*}
while by Theorem \ref{teor2}, we have
\begin{align*}
f(L_2,M_1)&-f(L_2,M_2)\nonumber\\[.2cm]
&=
\upint\!\!\!\iint\frac{f(x,y_1)-f(x,y_2)}{y_1-y_2}
\,d\E_{L_2}(x)\,d\E_{M_1}(y_1)(M_1-M_2)\,d\E_{M_2}(y_2).
\end{align*}
It remains to observe that
$$
f(L_1,M_1)-f(L_2,M_2)=(f(L_1,M_1)-f(L_2,M_1))+(f(L_2,M_1)-f(L_2,M_2)).
$$

Clearly, the set $\bigcup_{\s>0}(\E^\be_\s)_+(\R^2)$ is dense in the space $(\Bs)_+(\R^2)$.
It remains to observe that the set of functions $f$ in $(\Bs)_+(\R^2)$, for which inequality 
\rf{tryokhetazhka} holds is closed in the space $(\Bs)_+(\R^2)$.
This follows from Theorem \ref{lipschitseva_otsenka_dlya_Besova} and Theorem \ref{anp}. $\bl$

\medskip

Similarly, we can obtain the following result:

\begin{thm} 
\label{ba21}
Let $f\in(\Bs)_+(\R^2)$. Suppose that $L_1$, $L_2$, $M_1$ and $M_2$ are maximal dissipative operators such that 
$L_1-L_2\in\bS_2$ and $M_1-M_2\in\bS_2$. Then the following equality holds:
\begin{align*}
f(L_1,M_1)&-f(L_2,M_2)\nonumber\\[.2cm]
&=
\iint\!\!\upint\frac{f(x_1,y)-f(x_2,y)}{x_1-x_2}
\,d\E_{L_1}(x_1)(L_1-L_2)\,d\E_{L_2}(x_2)\,d\E_{M_2}(y),\nonumber\\[.2cm]
&+\upint\!\!\!\iint\frac{f(x,y_1)-f(x,y_2)}{y_1-y_2}
\,d\E_{L_1}(x)\,d\E_{M_1}(y_1)(M_1-M_2)\,d\E_{M_2}(y_2).
\end{align*}
\end{thm}

\

\

\
 
 \begin{footnotesize}
 
\noindent
\begin{tabular}{p{7cm}p{15cm}}
A.B. Aleksandrov & V.V. Peller \\
St.Petersburg Department & St.Petersburg State University\\
Steklov Institute of Mathematics  & Department of Mathematics and computer Science\\
Russian Academy of Sciences & Universitetskaya nab., 7/9\\
Fontanka 27, 191023 St.Petersburg &199034 St.Petersburg\\
Russia & Russia\\
email: aall54eexx@gmail.com\\
\\
&Department of Mathematics\\
&Michigan State University\\
&East Lansing, Michigan 48824\\
&USA\\
\\
&St.Petersburg Department\\
&Steklov Institute of Mathematics\\
&Russian Academy of Sciences\\
&Fontanka 27, 191023 St.Petersburg\\
&Russia\\
&and\\
&Peoples' Friendship University\\
&of Russia (RUDN University)\\
&6 Miklukho-Maklaya St., Moscow\\
& 117198, Russia\\
& email: peller@math.msu.edu
\end{tabular}

\end{footnotesize}

\end{document}